\documentclass[a4paper,12pt]{article}
\usepackage{latexsym} 
\usepackage{amssymb,amsmath} 
\usepackage{cite}

\newtheorem{definition}{Definition}

\begin{document}

\title
{
Local fractional Moisil-Teodorescu operator in quaternionic setting involving Cantor-type coordinate systems
}
\author{Juan Bory-Reyes$^{(1)}$ and Marco Antonio P\'{e}rez-de la Rosa$^{(2)}$}

\date{\small $^{(1)}$ ESIME-Zacatenco. Instituto Polit\'ecnico Nacional. CDMX. 07738. M\'exico.\\E-mail: juanboryreyes@yahoo.com\\
	$^{(2)}$ Department of Actuarial Sciences, Physics and Mathematics, Universidad de las Am\'{e}ricas Puebla.
	San Andr\'{e}s Cholula, Puebla. 72810. M\'{e}xico.\\ Email: marco.perez@udlap.mx}
\maketitle

\begin{abstract}
In this paper, general quaternionic structure are developed for the local fractional Moisil-Teodorescu operator in Cantor-type cylindrical and spherical coordinate systems. Two examples for the Helmholtz equation with local fractional derivatives on the Cantor sets are shown by making
use of this local fractional Moisil-Teodorescu operator.
\end{abstract}

\section{Introduction}
Vector analysis \cite{Tai, Ma}, a branch of mathematics that deals with scalar and vector quantities, has proven to be a powerful tool for developing a mathematical setting able to understanding several differential field equations in  different areas of engineering. 

One of the problems arising in vector analysis at Euclidean spaces, when dealing with operators such as divergence, gradient and curl for analysing the behaviour of scalar- and vector-valued functions, is that of understanding that each of these operators should fit all common differential operators and vectorial identities independently of the coordinate system considered.

It is well know that in vector analysis at $\mathbb{R}^3$ the vector product does not permit a formation of an algebra. In this sense, various attempts were made to construct such algebra. The vector algebra of quaternions, which marked the beginning of modern vector analysis, were discovered by William Rowan Hamilton, on the 16th of October $1843$. The fact that Hamilton started with triplets and ended up with quaternions implied immediately his concern about the special role of what he called vector: $\vec q= \sum_{j=1}^{3}q_{j}\mathbf{e}_{j}$, where $\{\mathbf{e}_{j}\}, j=1,2,3$ denotes standard orthonormal base of ${\mathbb R}^3$. 

The newly multiplication rules with the vectors of the base, as he used ever since:
\begin{equation}
\mathbf{e}^{2}_{1}=\mathbf{e}^{2}_{2}=\mathbf{e}^{2}_{3}=\mathbf{e}_{1}\mathbf{e}_{2}\mathbf{e}_{3}=-1
\end{equation}
contains the solution of the problem \cite{van, Lewis}. 

With the introduction of the Nabla-operator in Cartesian coordinates given by
\begin{equation}\label{nabla}
\nabla= \sum_{j=1}^{3}\mathbf{e}_{j}\partial_j,
\end{equation}
where $\partial_j$ is a partial derivative, Hamilton also invented the other essential technical ingredient for the vector calculus: the vector differential operator which is used to describe the gradient of a scalar function as well as the divergence and the curl of a vector valued
function.

There is a debate over whether it is correct that the action of the vector Laplacian is the action of the scalar Laplacian component by component, nevertheless this happens only for the case of the Cartesian coordinates and the result totally differ in any other orthogonal curvilinear system. \cite{Red, HiChi, Mo, Far}.

In the conventional approach, the Laplacian $\nabla^{2}$ is an operator that can operate on both scalar $q$ and vector $\vec q$ fields. The operator on a scalar can be written,
\begin{equation}
\nabla^{2} {q} = \nabla \cdot \nabla {q}
\end{equation}
which will produce another scalar field. Meanwhile, on a vector it can be expressed as
\begin{equation}
\nabla^{2} {\vec q} = {\nabla}(\nabla \cdot {\vec q}) - {\nabla} \times ({\nabla} \times {\vec q})
\end{equation}
which will produce another vector field. 

However, in Cartesian coordinates, both operators coincide with
\begin{equation}
\nabla^{2} = \sum_{j=1}^{3}\partial^{2}_j,
\end{equation}
where it is evident that operation on a scalar (vector) field transforms into a scalar (vector) field.

An overwhelming majority of physically meaningful problems may be simplified by use of non-Cartesian coordinates, especially orthogonal curvilinear coordinates.

Regarding the three-dimensional nature of many mathematical physics problems, to handle the techniques of quaternion analysis seems to be very perspective. Ideas and techniques from quaternionic analysis was originally introduced by Rudolf Fueter in \cite{Fu}, see also \cite{Fu1, Fu2}. So far, motivated by these works, mathematicians became interested in developing various approaches along classical lines in the study of this theory, see e.g. \cite{Ha, Dea, Sub}. 

The Moisil-Teodorescu operator (a determined first-order elliptic operator, that can be expressed in terms of the usual divergence, gradient and curl, see \cite{MT}) is nowadays considered to be a good analogue of the usual Cauchy-Riemann operator of complex analysis to the quaternionic setting and it is a square root of the scalar Laplace operator in $\mathbb{R}^3$. For a closely relation of  Moisil-Teodorescu operator to many mathematical models of important spatial physical phenomena we refer the reader to \cite{Zd, GuSp, Krav, KravShap}.

In \cite{Bory} the authors describe how the meaning of the Moisil-Teodorescu operator in open subsets of ${\mathbb R}^3$ involving orthogonal curvilinear coordinates may benefit from an approach within the framework of quaternionic analysis.

The enormous success and large application of local and non-Local fractional calculus (see \cite{Go, Gol, Rahmat, Yang, Yang2, Shi, Yang3, Hao, BaF, ORT, SZBCC}) is a strong motivation to show the quaternionic structure of local fractional Moisil-Teodorescu operator in Cantor-type cylindrical and spherical coordinates. This is the main goal of the paper. Additionally, the quaternionic Helmholtz equation associated with local fractional derivative operators involving the Cantor-type cylindrical and spherical coordinates is considered.
\section{Preliminaries}

\subsection{Few aspects of local fractional calculus}

We begin our study by recalling some facts on local fractional calculus that can be found in \cite{Yang,Shi}. 

\begin{definition}
Let $f(x)$ be a function defined on a fractal set of fractal dimension $\alpha$ $(0 < \alpha < 1)$, the function $f(x)$ is said to be local fractional continuous at $x=x_0$ if for each $\epsilon >0$, there exists a corresponding $\delta >0$ such that
\[|f(x)-f(x_0)|<\epsilon^{\alpha},\]	
whenever $0<|x-x_0|<\delta$.
\end{definition}

The set of all local fractional continuous functions on the interval $(a, b)$ will be denoted by $C_{\alpha}(a,b)$.

\begin{definition}
	The local fractional derivative of $f(x)\in C_{\alpha}(a,b)$ of order $\alpha$ $(0<\alpha<1)$ at $x=x_0$ is defined as
	\[D^{(\alpha)}f(x_0)=f^{(\alpha)}(x_0)=\left.\frac{d^{\alpha}f(x)}{dx^{\alpha}}\right|_{x=x_0}=\lim_{x\to x_0}\frac{\Delta^{\alpha}\left[f(x)-f(x_0)\right]}{(x-x_0)^{\alpha}},\]
	provided the limit exists, where $\displaystyle \Delta^{\alpha}\left[f(x)-f(x_0)\right]\cong \Gamma(1+\alpha)\left[f(x)-f(x_0)\right]$ with Euler's gamma function $\Gamma(1+\alpha):=\int_{0}^{\infty}\mu^{\alpha-1}\,\exp(-\mu)\,d\mu$.
\end{definition}

The set of all functions such that its local derivative exists for all $x$ on the interval $(a, b)$ will be denoted by $D_{\alpha}(a,b)$. The set $D_{\alpha}(a,b)$ is called a $\alpha$-local fractional derivative set.

Let $f,g\in D_{\alpha}(a,b)$, then
\begin{align}
&D^{(\alpha)}\left[f(x)\pm g(x)\right]=D^{(\alpha)}f(x)\pm D^{(\alpha)}g(x),\\
&D^{(\alpha)}\left[f(x)g(x)\right]=g(x)D^{(\alpha)}f(x)+f(x)D^{(\alpha)}g(x),\\
&D^{(\alpha)}\left[\frac{f(x)}{g(x)}\right]=\frac{\left[D^{(\alpha)}f(x)\right]g(x)-f(x)\left[D^{(\alpha)}f(x)\right]}{g^2(x)},
\end{align}
provided $g(x)\neq0$.

Suppose that $f(x)=(\phi\circ\varphi)(x)$, $x\in(a,b)$. Then, we have
\[f^{(\alpha)}(x)=\phi^{(\alpha)}(\varphi(x))\left[\varphi^{(1)}(x)\right]^{\alpha},\]
provided $\phi^{(\alpha)}(\varphi(x))$ and $\varphi^{(1)}(x)$ exist.

The generalized functions defined on Cantor sets are given by
\begin{equation}
E_{\alpha}\left(x^{\alpha}\right)=\sum_{k=0}^{\infty}\frac{x^{k\alpha}}{\Gamma(1+k\alpha)}.
\end{equation}

The sine and cosine functions on a fractal set are given, respectively, by
\begin{align}
\sin_{\alpha}\left(x^{\alpha}\right)&=\frac{E_{\alpha}\left(i^{\alpha}x^{\alpha}\right)-E_{\alpha}\left(-i^{\alpha}x^{\alpha}\right)}{2i^{\alpha}}=\sum_{k=0}^{\infty}\frac{(-1)^kx^{(2k+1)\alpha}}{\Gamma(1+(2k+1)\alpha)},\\
\cos_{\alpha}\left(x^{\alpha}\right)&=\frac{E_{\alpha}\left(i^{\alpha}x^{\alpha}\right)+E_{\alpha}\left(-i^{\alpha}x^{\alpha}\right)}{2}=\sum_{k=0}^{\infty}\frac{(-1)^kx^{2k\alpha}}{\Gamma(1+2k\alpha)},
\end{align}
where $x\in\mathbb{R}$, $0<\alpha<1$, and $i^{\alpha}$ is a imaginary unit of a fractal set.

The previously defined functions satisfy the following identities:
\begin{align}
&D^{(\alpha)}\left[\frac{x^{n\alpha}}{\Gamma(1+n\alpha)}\right]=\frac{x^{(n-1)\alpha}}{\Gamma(1+(n-1)\alpha)},\\
&D^{(\alpha)}\left[E_{\alpha}\left(x^{\alpha}\right)\right]=E_{\alpha}\left(x^{\alpha}\right),\\
%&D^{(\alpha)}\left[E_{\alpha}\left(cx^{\alpha}\right)\right]=cE_{\alpha}\left(cx^{\alpha}\right),\\
%&D^{(\alpha)}\left[E_{\alpha}\left(x^{2\alpha}\right)\right]=(2x)^{\alpha}E_{\alpha}\left(x^{2\alpha}\right),\\
&D^{(\alpha)}\left[\sin_{\alpha}\left(x^{\alpha}\right)\right]=\cos_{\alpha}\left(x^{\alpha}\right),\\
&D^{(\alpha)}\left[\cos_{\alpha}\left(x^{\alpha}\right)\right]=-\sin_{\alpha}\left(x^{\alpha}\right).
%&D^{(\alpha)}\left[\sin_{\alpha}^{\alpha}\left(x^{\alpha}\right)\right]=\Gamma(1+\alpha)\cos_{\alpha}\left(x^{\alpha}\right),\\
%&D^{(\alpha)}\left[\cos_{\alpha}^{\alpha}\left(x^{\alpha}\right)\right]=-\Gamma(1+\alpha)\sin_{\alpha}\left(x^{\alpha}\right).
\end{align}

\subsection{Rudiments of quaternionic analysis}

In this subsection, we follow \cite{KravShap} in reviewing some standard facts on quaternionic analysis to be used in this paper. 

We work with the skew-field (a complex non-commutative, associative algebra with zero divisors)  $\mathbb H(\mathbb C)$ of complex quaternions, i.e., each $a\in \mathbb H(\mathbb C)$ is of the form $q = \sum_{k=0}^3 q_k\mathbf{i}_k$, with $\{q_k\}\subset\mathbb C$; $\mathbf{i}_0 = 1$ and $\mathbf{i}_1, \mathbf{i}_2, \mathbf{i}_3$ stand for the quaternionic imaginary units. By definition, the complex imaginary unit in $\mathbb{C}$, denotes by $i$, commutes with all the quaternionic imaginary units.

For $q = \sum_{k=0}^3 q_k\mathbf{i}_k\in\mathbb H(\mathbb C)$, we will write $q_0 =: \text{Sc}(q)$, $\vec{q} := \sum_{k=1}^3 q_k\mathbf{i}_k =: \text{Vec}(q)$, so $q = q_0 + \vec{q}$. We call $q_0$ and $\vec{q}$ the scalar and vector parts of $q$ respectively. Then $\{ \text{Vec}(q) : q\in\mathbb H(\mathbb C)\}$ is identified with $\mathbb C^3$.

For any $p,q\in\mathbb H(\mathbb C)$:
\[
p\,q := p_0\,q_0 - \langle\vec p,\vec q\rangle + p_0\,\vec q + q_0\,\vec p + [\vec p,\vec q],
\]
where
\[
\langle\vec p,\vec q\rangle := \sum_{k=1}^3 p_k\,q_k,\,\, [\vec p,\vec q]:=
\left |
\begin{array}{rrr}
\mathbf{i}_1 & \mathbf{i}_2 & \mathbf{i}_3\\
p_1& p_2 & p_3\\
q_1& q_2 & q_3
\end{array}
\right |.
\]
In particular, if $p_0=q_0=0$ then $p\,q :=  - \langle\vec p,\vec q\rangle + [\vec p,\vec q]$.

The Moisil-Teodorescu operator, denoted by $D_{MT}$, is defined to be:
\begin{equation}
D_{MT}[f]:=\mathbf{i}_1\frac{\partial f}{\partial x}+\mathbf{i}_2\frac{\partial f}{\partial y}+\mathbf{i}_3\frac{\partial f}{\partial z}.
\end{equation}

It is worth noting that $D_{MT}$ be resemblance to $\nabla$ and factorizes the scalar Laplacian. In fact, it holds
\begin{equation}
-D_{MT}^2=\Delta_{\mathbb{R}^4},
\end{equation}
which implies several advantages in the applications to physical problems.

The operator $\Delta_{\mathbb{R}^4}$ is a scalar operator, its acts separately on every coordinate function $f_k$ of $f$ as
$$\Delta_{\mathbb{R}^4}[f]:=\Delta[f_0]+\mathbf{i}_1\Delta[f_1]+\mathbf{i}_2\Delta[f_2]+\mathbf{i}_3\Delta[f_3].$$

This property guarantees that any hyperholomorphic function is also harmonic.

For an $\mathbb{H}(\mathbb{C})$-valued function $f:=f_0+ \vec{f}$ the action of the operator $D_{MT}$ can be
represented as follows
\begin{equation}
D_{MT}[f]=-\mathrm{div}[\vec{f}]+\mathrm{grad}[f_0]+\mathrm{curl}[\vec{f}].
\end{equation}
This is an immediate consequence of the quaternionic product.

Moisil-Teodorescu operator can act on the right in which case the notation $D_{MT}^r[f]$  for the same $f=f_0+\vec f$ means that:
\begin{equation}
D_{MT}^r[f]:=\frac{\partial f}{\partial x}\mathbf{i}_1+\frac{\partial f}{\partial y}\mathbf{i}_2+\frac{\partial f}{\partial z}\mathbf{i}_3.
\end{equation}

\section{Local fractional Moisil-Teodorescu operator in Cantor-type coordinates}

\subsection{The Cantor-type cylindrical coordinate system}

Consider the coordinate system of the Cantor-type cylindrical coordinates developed in \cite{Hao,Yang2,Yang,Shi,Yang3}.

Consider the Cantor-type cylindrical coordinates
\begin{equation}
	\begin{cases}
		x^{\alpha}=r^{\alpha}\cos_{\alpha}\left(\theta^{\alpha}\right),\\
		y^{\alpha}=r^{\alpha}\sin_{\alpha}\left(\theta^{\alpha}\right),\\
		z^{\alpha}=z^{\alpha}.
	\end{cases}
\end{equation}

If we denote by
\begin{equation}
	\begin{cases}
		\mathbf{e}_{r}^{\alpha}=\cos_{\alpha}\left(\theta^{\alpha}\right)\mathbf{i}_{1}^{\alpha}+\sin_{\alpha}\left(\theta^{\alpha}\right)\mathbf{i}_{2}^{\alpha},\\
		\mathbf{e}_{\theta}^{\alpha}=-\sin_{\alpha}\left(\theta^{\alpha}\right)\mathbf{i}_{1}^{\alpha}+\cos_{\alpha}\left(\theta^{\alpha}\right)\mathbf{i}_{2}^{\alpha},\\
		\mathbf{e}_{\psi}^{\alpha}=\mathbf{i}_{3}^{\alpha},
	\end{cases}
\end{equation}
then the vectorial operations gradient, divergence and curl in orthogonal curvilinear coordinates can be expressed as follows:

If $f=f_0$ a scalar function, then one has
\begin{equation}
	\mathrm{grad}^{(\alpha)}[f_0]=\frac{\partial^{\alpha} f_0}{\partial r^{\alpha}}\,\mathbf{e}_{r}^{\alpha}+\frac{1}{r^{\alpha}}\frac{\partial^{\alpha} f_0}{\partial \theta^{\alpha}}\,\mathbf{e}_{\theta}^{\alpha}+\frac{\partial^{\alpha} f_0}{\partial z^{\alpha}}\,\mathbf{e}_{z}^{\alpha}.
\end{equation}

If $f=\vec{f}=f_1\mathbf{e}_{r}^{\alpha}+f_2\mathbf{e}_{\theta}^{\alpha}+f_3\mathbf{e}_{z}^{\alpha}$, then we have

\begin{equation}
	\mathrm{div}^{(\alpha)}[\vec{f}]=\frac{\partial^{\alpha} f_1}{\partial r^{\alpha}}+\frac{1}{r^{\alpha}}\frac{\partial^{\alpha} f_2}{\partial \theta^{\alpha}}+\frac{f_1}{r^{\alpha}}+\frac{\partial^{\alpha} f_3}{\partial z^{\alpha}}
\end{equation}

and

\begin{align}
	\mathrm{curl}^{(\alpha)}[\vec{f}]=&\left[\frac{1}{r^{\alpha}}\frac{\partial^{\alpha} f_3}{\partial \theta^{\alpha}}-\frac{\partial^{\alpha} f_2}{\partial z^{\alpha}}\right]\,\mathbf{e}_{r}^{\alpha}+\left[\frac{\partial^{\alpha} f_1}{\partial z^{\alpha}}-\frac{\partial^{\alpha} f_3}{\partial r^{\alpha}}\right]\,\mathbf{e}_{\theta}^{\alpha}\notag\\
	&+\left[\frac{\partial^{\alpha} f_2}{\partial r^{\alpha}}-\frac{1}{r^{\alpha}}\frac{\partial^{\alpha} f_1}{\partial \theta^{\alpha}}+\frac{f_2}{r^{\alpha}}\right]\,\mathbf{e}_{z}^{\alpha}.
\end{align}

\subsection{The Cantor-type spherical coordinate system}

Consider the coordinate system of the Cantor-type spherical coordinates developed in \cite{Rahmat, Shi, Yang, Yang3}.

Consider the Cantor-type spherical coordinates
\begin{equation}
	\begin{cases}
		x^{\alpha}=r^{\alpha}\sin_{\alpha}\left(\theta^{\alpha}\right)\cos_{\alpha}\left(\psi^{\alpha}\right),\\
		y^{\alpha}=r^{\alpha}\sin_{\alpha}\left(\theta^{\alpha}\right)\sin_{\alpha}\left(\psi^{\alpha}\right),\\
		z^{\alpha}=r^{\alpha}\cos_{\alpha}\left(\theta^{\alpha}\right).
	\end{cases}
\end{equation}

If we denote by
\begin{equation}
	\begin{cases}
		\mathbf{e}_{r}^{\alpha}=\sin_{\alpha}\left(\theta^{\alpha}\right)\cos_{\alpha}\left(\psi^{\alpha}\right)\mathbf{i}_{1}^{\alpha}+\sin_{\alpha}\left(\theta^{\alpha}\right)\sin_{\alpha}\left(\psi^{\alpha}\right)\mathbf{i}_{2}^{\alpha}+\cos_{\alpha}\left(\theta^{\alpha}\right)\mathbf{i}_{3}^{\alpha},\\
		\mathbf{e}_{\theta}^{\alpha}=\cos_{\alpha}\left(\theta^{\alpha}\right)\cos_{\alpha}\left(\psi^{\alpha}\right)\mathbf{i}_{1}^{\alpha}+\cos_{\alpha}\left(\theta^{\alpha}\right)\sin_{\alpha}\left(\psi^{\alpha}\right)\mathbf{i}_{2}^{\alpha}-\sin_{\alpha}\left(\theta^{\alpha}\right)\mathbf{i}_{3}^{\alpha},\\
		\mathbf{e}_{\psi}^{\alpha}=-\sin_{\alpha}\left(\psi^{\alpha}\right)\mathbf{i}_{1}^{\alpha}+\cos_{\alpha}\left(\psi^{\alpha}\right)\mathbf{i}_{2}^{\alpha},
	\end{cases}
\end{equation}
then the vectorial operations gradient, divergence and curl in orthogonal curvilinear coordinates can be expressed as follows:

If $f=f_0$ a scalar function, then one has
\begin{equation}
	\mathrm{grad}^{(\alpha)}[f_0]=\frac{\partial^{\alpha} f_0}{\partial r^{\alpha}}\,\mathbf{e}_{r}^{\alpha}+\frac{1}{r^{\alpha}}\frac{\partial^{\alpha} f_0}{\partial \theta^{\alpha}}\,\mathbf{e}_{\theta}^{\alpha}+\frac{1}{r^{\alpha}\sin_{\alpha}\left(\theta^{\alpha}\right)}\frac{\partial^{\alpha} f_0}{\partial \psi^{\alpha}}\,\mathbf{e}_{\psi}^{\alpha}.
\end{equation}

If $f=\vec{f}=f_1\mathbf{e}_{r}^{\alpha}+f_2\mathbf{e}_{\theta}^{\alpha}+f_3\mathbf{e}_{\psi}^{\alpha}$, then we have

\begin{equation}
	\mathrm{div}^{(\alpha)}[\vec{f}]=\frac{\partial^{\alpha} f_1}{\partial r^{\alpha}}+\frac{2f_1}{r^{\alpha}}+\frac{1}{r^{\alpha}}\frac{\partial^{\alpha} f_2}{\partial \theta^{\alpha}}+\frac{1}{r^{\alpha}\sin_{\alpha}\left(\theta^{\alpha}\right)}\left[\frac{\partial^{\alpha} f_3}{\partial \psi^{\alpha}}+f_2\cos_{\alpha}\left(\theta^{\alpha}\right)\right],
\end{equation}

and

\begin{align}
	\mathrm{curl}^{(\alpha)}[\vec{f}]=&\left[\frac{1}{r^{\alpha}}\frac{\partial^{\alpha} f_3}{\partial \theta^{\alpha}}-\frac{1}{r^{\alpha}\sin_{\alpha}\left(\theta^{\alpha}\right)}\frac{\partial^{\alpha} f_2}{\partial \psi^{\alpha}}+\frac{f_3\cos_{\alpha}\left(\theta^{\alpha}\right)}{r^{\alpha}\sin_{\alpha}\left(\theta^{\alpha}\right)}\right]\,\mathbf{e}_{r}^{\alpha}\notag\\
	&+\left[\frac{1}{r^{\alpha}\sin_{\alpha}\left(\theta^{\alpha}\right)}\frac{\partial^{\alpha} f_1}{\partial \psi^{\alpha}}-\frac{\partial^{\alpha} f_3}{\partial r^{\alpha}}-\frac{f_3}{r^{\alpha}}\right]\,\mathbf{e}_{\theta}^{\alpha}\notag\\
	&+\left[\frac{\partial^{\alpha} f_2}{\partial r^{\alpha}}-\frac{1}{r^{\alpha}}\frac{\partial^{\alpha} f_1}{\partial \theta^{\alpha}}+\frac{f_2}{r^{\alpha}}\right]\,\mathbf{e}_{\psi}^{\alpha}.
\end{align}

\subsection{Local fractional Moisil-Teodorescu operator in\\ Can\-tor-type cylindrical and spherical coordinates}

Using the results from the previous two sections we have that the local fractional Moisil-Teodorescu operator for a quaternionic-valued function in Cantor-type cylindrical coordinates is given by:
\begin{align}
D_{MT}^{(\alpha)}[f]=&-\mathrm{div}^{(\alpha)}[\vec{f}]+\mathrm{grad}^{(\alpha)}[f_0]+\mathrm{curl}^{(\alpha)}[\vec{f}]\notag\\
=&-\left[\frac{\partial^{\alpha} f_1}{\partial r^{\alpha}}+\frac{1}{r^{\alpha}}\frac{\partial^{\alpha} f_2}{\partial \theta^{\alpha}}+\frac{f_1}{r^{\alpha}}+\frac{\partial^{\alpha} f_3}{\partial z^{\alpha}}\right]\notag\\
&+\left[\frac{\partial^{\alpha} f_0}{\partial r^{\alpha}}+\frac{1}{r^{\alpha}}\frac{\partial^{\alpha} f_3}{\partial \theta^{\alpha}}-\frac{\partial^{\alpha} f_2}{\partial z^{\alpha}}\right]\,\mathbf{e}_{r}^{\alpha}\notag\\
&+\left[\frac{1}{r^{\alpha}}\frac{\partial^{\alpha} f_0}{\partial \theta^{\alpha}}+\frac{\partial^{\alpha} f_1}{\partial z^{\alpha}}-\frac{\partial^{\alpha} f_3}{\partial r^{\alpha}}\right]\,\mathbf{e}_{\theta}^{\alpha}\notag\\
&+\left[\frac{\partial^{\alpha} f_0}{\partial z^{\alpha}}+\frac{\partial^{\alpha} f_2}{\partial r^{\alpha}}-\frac{1}{r^{\alpha}}\frac{\partial^{\alpha} f_1}{\partial \theta^{\alpha}}+\frac{f_2}{r^{\alpha}}\right]\,\mathbf{e}_{z}^{\alpha},
\end{align}
while in Cantor-type spherical coordinates is given by:
\begin{align}
	D_{MT}^{(\alpha)}[f]=&-\mathrm{div}^{(\alpha)}[\vec{f}]+\mathrm{grad}^{(\alpha)}[f_0]+\mathrm{curl}^{(\alpha)}[\vec{f}]\notag\\
	=&-\left[\frac{\partial^{\alpha} f_1}{\partial r^{\alpha}}+\frac{2f_1}{r^{\alpha}}+\frac{1}{r^{\alpha}}\frac{\partial^{\alpha} f_2}{\partial \theta^{\alpha}}+\frac{1}{r^{\alpha}\sin_{\alpha}\left(\theta^{\alpha}\right)}\left(\frac{\partial^{\alpha} f_3}{\partial \psi^{\alpha}}+f_2\cos_{\alpha}\left(\theta^{\alpha}\right)\right)\right]\notag\\
	&+\left[\frac{\partial^{\alpha} f_0}{\partial r^{\alpha}}+\frac{1}{r^{\alpha}}\frac{\partial^{\alpha} f_3}{\partial \theta^{\alpha}}-\frac{1}{r^{\alpha}\sin_{\alpha}\left(\theta^{\alpha}\right)}\frac{\partial^{\alpha} f_2}{\partial \psi^{\alpha}}+\frac{f_3\cos_{\alpha}\left(\theta^{\alpha}\right)}{r^{\alpha}\sin_{\alpha}\left(\theta^{\alpha}\right)}\right]\,\mathbf{e}_{r}^{\alpha}\notag\\
	&+\left[\frac{1}{r^{\alpha}}\frac{\partial^{\alpha} f_0}{\partial \theta^{\alpha}}+\frac{1}{r^{\alpha}\sin_{\alpha}\left(\theta^{\alpha}\right)}\frac{\partial^{\alpha} f_1}{\partial \psi^{\alpha}}-\frac{\partial^{\alpha} f_3}{\partial r^{\alpha}}-\frac{f_3}{r^{\alpha}}\right]\,\mathbf{e}_{\theta}^{\alpha}\notag\\
	&+\left[\frac{1}{r^{\alpha}\sin_{\alpha}\left(\theta^{\alpha}\right)}\frac{\partial^{\alpha} f_0}{\partial \psi^{\alpha}}+\frac{\partial^{\alpha} f_2}{\partial r^{\alpha}}-\frac{1}{r^{\alpha}}\frac{\partial^{\alpha} f_1}{\partial \theta^{\alpha}}+\frac{f_2}{r^{\alpha}}\right]\,\mathbf{e}_{\psi}^{\alpha}.
\end{align}

\subsection{Local fractional Laplacian in Cantor-type cylindrical and spherical coordinates}

The local fractional Moisil-Teodorescu operator $D_{MT}^{(\alpha)}$ factorizes the Laplacian $\Delta_{\mathbb H}^{(\alpha)}$ as follows:
\begin{equation}
	\left(D_{MT}^{(\alpha)}\right)^2[f]=-\Delta_{\mathbb H}^{(\alpha)}[f],
\end{equation}
where
\begin{equation}
\Delta_{\mathbb H}^{(\alpha)}[f]:=\Delta_0^{(\alpha)}[f_0]+\vec{\Delta}^{(\alpha)}\left[\vec{f}\right],
\end{equation}
with
$
\Delta_0^{(\alpha)}[f_0]:=\mathrm{div}^{(\alpha)}\left[\mathrm{grad}^{(\alpha)}[f_0]\right]$ and $\vec{\Delta}^{(\alpha)}\left[\vec{f}\right]:=\mathrm{grad}^{(\alpha)}\left[\mathrm{div}^{(\alpha)}\left[\vec{f}\right]\right]-\mathrm{curl}^{(\alpha)}\left[\mathrm{curl}^{(\alpha)}\left[\vec{f}\right]\right]$.

Then, in Cantor-type cylindrical coordinates the previous operators take the form:
\begin{equation}\label{scalar Laplacian cylindrical}
\Delta_0^{(\alpha)}[f_0]=\frac{\partial^{2\alpha} f_0}{\partial r^{2\alpha}}+\frac{1}{r^{2\alpha}}\frac{\partial^{2\alpha} f_0}{\partial \theta^{2\alpha}}+\frac{1}{r^{\alpha}}\frac{\partial^{\alpha} f_0}{\partial r^{\alpha}}+\frac{\partial^{2\alpha} f_0}{\partial z^{2\alpha}}
\end{equation}
and
\begin{align}
	\vec{\Delta}^{(\alpha)}\left[\vec{f}\right]&=\left[\frac{\partial^{2\alpha} f_1}{\partial r^{2\alpha}}+\frac{1}{r^{\alpha}}\frac{\partial^{\alpha} f_1}{\partial r^{\alpha}}-\frac{f_1}{r^{2\alpha}}+\frac{1}{r^{2\alpha}}\frac{\partial^{2\alpha} f_1}{\partial \theta^{2\alpha}}-\frac{2}{r^{2\alpha}}\frac{\partial^{\alpha} f_2}{\partial \theta^{\alpha}}+\frac{\partial^{2\alpha} f_1}{\partial z^{2\alpha}}\right]\,\mathbf{e}_{r}^{\alpha}\notag\\
	&+\left[\frac{\partial^{2\alpha} f_2}{\partial r^{2\alpha}}+\frac{1}{r^{\alpha}}\frac{\partial^{\alpha} f_2}{\partial r^{\alpha}}-\frac{f_2}{r^{2\alpha}}+\frac{1}{r^{2\alpha}}\frac{\partial^{2\alpha} f_2}{\partial \theta^{2\alpha}}+\frac{2}{r^{2\alpha}}\frac{\partial^{\alpha} f_1}{\partial \theta^{\alpha}}+\frac{\partial^{2\alpha} f_2}{\partial z^{2\alpha}}\right]\,\mathbf{e}_{\theta}^{\alpha}\notag\\
	&+\left[\frac{\partial^{2\alpha} f_3}{\partial r^{2\alpha}}+\frac{1}{r^{\alpha}}\frac{\partial^{\alpha} f_3}{\partial r^{\alpha}}+\frac{1}{r^{2\alpha}}\frac{\partial^{2\alpha} f_3}{\partial \theta^{2\alpha}}+\frac{\partial^{2\alpha} f_3}{\partial z^{\alpha}}\right]\,\mathbf{e}_{z}^{\alpha}\label{vector Laplacian cylindrical}
\end{align}

Moreover, in Cantor-type spherical coordinates these operators take the form:
\begin{align}
	\Delta_0^{(\alpha)}[f_0]&=\frac{\partial^{2\alpha} f_0}{\partial r^{2\alpha}}+\frac{2}{r^{\alpha}}\frac{\partial^{\alpha} f_0}{\partial r^{\alpha}}+\frac{1}{r^{2\alpha}}\frac{\partial^{2\alpha} f_0}{\partial \theta^{2\alpha}}+\frac{\cos_{\alpha}\left(\theta^{\alpha}\right)}{r^{2\alpha}\sin_{\alpha}\left(\theta^{\alpha}\right)}\frac{\partial^{\alpha} f_0}{\partial \theta^{\alpha}}\notag\\
	&+\frac{1}{r^{2\alpha}\sin_{\alpha}^{2\alpha}\left(\theta^{\alpha}\right)}\frac{\partial^{2\alpha} f_0}{\partial \psi^{2\alpha}},\label{scalar Laplacian spherical}
\end{align}
and
\begin{align}
	&\vec{\Delta}^{(\alpha)}\left[\vec{f}\right]=\left[\frac{\partial^{2\alpha} f_1}{\partial r^{2\alpha}}+\frac{2}{r^{\alpha}}\frac{\partial^{\alpha} f_1}{\partial r^{\alpha}}-\frac{2}{r^{2\alpha}}f_1+\frac{1}{r^{2\alpha}}\frac{\partial^{2\alpha} f_1}{\partial \theta^{2\alpha}}+\frac{\cos_{\alpha}\left(\theta^{\alpha}\right)}{r^{2\alpha}\sin_{\alpha}\left(\theta^{\alpha}\right)}\frac{\partial^{\alpha} f_1}{\partial \theta^{\alpha}}\right.\notag\\
	&+\left.\frac{1}{r^{2\alpha}\sin_{\alpha}^{2\alpha}\left(\theta^{\alpha}\right)}\frac{\partial^{2\alpha} f_1}{\partial \psi^{2\alpha}}-\frac{2}{r^{2\alpha}}\frac{\partial^{\alpha} f_2}{\partial \theta^{\alpha}}-\frac{2\cos_{\alpha}\left(\theta^{\alpha}\right)}{r^{2\alpha}\sin_{\alpha}\left(\theta^{\alpha}\right)}f_2-\frac{2}{r^{2\alpha}\sin_{\alpha}\left(\theta^{\alpha}\right)}\frac{\partial^{\alpha} f_3}{\partial \psi^{\alpha}}\right]\,\mathbf{e}_{r}^{\alpha}\notag\\
	&+\left[\frac{\partial^{2\alpha} f_2}{\partial r^{2\alpha}}+\frac{2}{r^{\alpha}}\frac{\partial^{\alpha} f_2}{\partial r^{\alpha}}-\frac{1}{r^{2\alpha}\sin_{\alpha}^{2\alpha}\left(\theta^{\alpha}\right)}f_2+\frac{1}{r^{2\alpha}}\frac{\partial^{2\alpha} f_2}{\partial \theta^{2\alpha}}+\frac{\cos_{\alpha}\left(\theta^{\alpha}\right)}{r^{2\alpha}\sin_{\alpha}\left(\theta^{\alpha}\right)}\frac{\partial^{\alpha} f_2}{\partial \theta^{\alpha}}\right.\notag\\
	&+\left.\frac{1}{r^{2\alpha}\sin_{\alpha}^{2\alpha}\left(\theta^{\alpha}\right)}\frac{\partial^{2\alpha} f_2}{\partial \psi^{2\alpha}}+\frac{2}{r^{2\alpha}}\frac{\partial^{\alpha} f_1}{\partial \theta^{\alpha}}-\frac{2\cos_{\alpha}\left(\theta^{\alpha}\right)}{r^{2\alpha}\sin_{\alpha}^{2\alpha}\left(\theta^{\alpha}\right)}\frac{\partial^{\alpha} f_3}{\partial \psi^{\alpha}}\right]\,\mathbf{e}_{\theta}^{\alpha}\notag\\
	&+\left[\frac{\partial^{2\alpha} f_3}{\partial r^{2\alpha}}+\frac{2}{r^{\alpha}}\frac{\partial^{\alpha} f_3}{\partial r^{\alpha}}-\frac{1}{r^{2\alpha}\sin_{\alpha}^{2\alpha}\left(\theta^{\alpha}\right)}f_3+\frac{1}{r^{2\alpha}}\frac{\partial^{2\alpha} f_3}{\partial \theta^{2\alpha}}+\frac{\cos_{\alpha}\left(\theta^{\alpha}\right)}{r^{2\alpha}\sin_{\alpha}\left(\theta^{\alpha}\right)}\frac{\partial^{\alpha} f_3}{\partial \theta^{\alpha}}\right.\notag\\
	&+\left.\frac{1}{r^{2\alpha}\sin_{\alpha}^{2\alpha}\left(\theta^{\alpha}\right)}\frac{\partial^{2\alpha} f_3}{\partial \psi^{2\alpha}}+\frac{2}{r^{2\alpha}\sin_{\alpha}\left(\theta^{\alpha}\right)}\frac{\partial^{\alpha} f_1}{\partial \psi^{\alpha}}+\frac{2\cos_{\alpha}\left(\theta^{\alpha}\right)}{r^{2\alpha}\sin_{\alpha}^{2\alpha}\left(\theta^{\alpha}\right)}\frac{\partial^{\alpha} f_2}{\partial \psi^{\alpha}}\right]\,\mathbf{e}_{\psi}^{\alpha}.\label{vector Laplacian spherical}
\end{align}

\subsection{Local fractional Bitsadze operator}

The local fractional quaternionic Bitsadze operator reads as 
\begin{equation}
\widetilde{\Delta_{\mathbb H}^{(\alpha)}}[f]:=\Delta_0^{(\alpha)}[f_0]+\widetilde{\vec{\Delta}}^{(\alpha)}\left[\vec{f}\right],
\end{equation}
with $
\widetilde{\vec{\Delta}}^{(\alpha)}\left[\vec{f}\right]:=\mathrm{grad}^{(\alpha)}\left[\mathrm{div}^{(\alpha)}\left[\vec{f}\right]\right]+\mathrm{curl}^{(\alpha)}\left[\mathrm{curl}^{(\alpha)}\left[\vec{f}\right]\right]$.

Notice that $D_{MT}^{(\alpha)}$ and $D_{MT}^{r(\alpha)}$ factorize the local fractional quaternionic Bitsadze operator $\widetilde{\Delta_{\mathbb H}^{(\alpha)}}$ as follow:
\begin{equation}
	D_{MT}^{(\alpha)} D_{MT}^{r(\alpha)}[f]=-\widetilde{\Delta_{\mathbb H}^{(\alpha)}}[f].
\end{equation}

Our purpose is to derive the local fractional quaternionic Bitsadze operator on the Cantor sets by using the Cantor-type cylindrical and spherical coordinates, which extend the quaternionic Bitsadze operator of \cite{Bory} based upon the standard derivative operators. 

We are reduced to handle the corresponding local fractional expressions of $\widetilde{\vec{\Delta}}^{(\alpha)}$ in Cantor-type cylindrical and spherical coordinates

\subsubsection*{Cantor-type cylindrical  coordinates}

\begin{align}
\widetilde{\vec{\Delta}}^{(\alpha)}\left[\vec{f}\right]&=\left[\frac{\partial^{2\alpha} f_1}{\partial r^{2\alpha}}+\frac{2}{r^{\alpha}}\frac{\partial^{2\alpha} f_2}{\partial r^{\alpha}\partial\theta^{\alpha}}+\frac{1}{r^{\alpha}}\frac{\partial^{\alpha} f_1}{\partial r^{\alpha}}-\frac{f_1}{r^{2\alpha}}+2\frac{\partial^{2\alpha} f_3}{\partial r^{\alpha} \partial z^{\alpha}}\right.\notag\\
&\left.-\frac{1}{r^{2\alpha}}\frac{\partial^{2\alpha} f_1}{\partial  \theta^{2\alpha}}-\frac{\partial^{2\alpha} f_1}{\partial z^{2\alpha} }\right]\,\mathbf{e}_{r}^{\alpha}\notag\\
&+\left[\frac{2}{r^{\alpha}}\frac{\partial^{2\alpha} f_1}{\partial r^{\alpha}\partial\theta^{\alpha}}+\frac{1}{r^{2\alpha}}\frac{\partial^{2\alpha} f_2}{\partial \theta^{2\alpha}}+\frac{2}{r^{\alpha}}\frac{\partial^{2\alpha} f_3}{\partial \theta^{\alpha}\partial z^{\alpha}}-\frac{\partial^{2\alpha} f_2}{\partial z^{2\alpha}}-\frac{\partial^{2\alpha} f_2}{\partial r^{2\alpha}}\right.\notag\\
&-\left.\frac{1}{r^{\alpha}}\frac{\partial^{\alpha} f_2}{\partial r^{\alpha}}+\frac{f_2}{r^{2\alpha}}\right]\,\mathbf{e}_{\theta}^{\alpha}\notag\\
&+\left[2\frac{\partial^{2\alpha} f_1}{\partial r^{\alpha}\partial z^{\alpha}}+\frac{2}{r^{\alpha}}\frac{\partial^{2\alpha} f_2}{\partial \theta^{\alpha}\partial z^{\alpha}}+\frac{2}{r^{\alpha}}\frac{\partial^{\alpha} f_1}{\partial z^{\alpha}}+\frac{\partial^{2\alpha} f_3}{\partial z^{2\alpha}}-\frac{\partial^{2\alpha} f_3}{\partial r^{2\alpha}}\right.\notag\\
&-\left.\frac{1}{r^{2\alpha}}\frac{\partial^{2\alpha} f_3}{\partial \theta^{2\alpha}}-\frac{1}{r^{\alpha}}\frac{\partial^{\alpha} f_3}{\partial r^{\alpha}}\right]\,\mathbf{e}_{z}^{\alpha}.
\end{align}

\subsubsection*{Cantor-type spherical coordinates}

\begin{align}
	\widetilde{\vec{\Delta}}^{(\alpha)}\left[\vec{f}\right]&=\left[\frac{\partial^{2\alpha} f_1}{\partial r^{2\alpha}}+\frac{2}{r^{\alpha}}\frac{\partial^{\alpha} f_1}{\partial r^{\alpha}}-\frac{2}{r^{2\alpha}}f_1-\frac{1}{r^{2\alpha}}\frac{\partial^{2\alpha} f_1}{\partial \theta^{2\alpha}}-\frac{\cos_{\alpha}\left(\theta^{\alpha}\right)}{r^{2\alpha}\sin_{\alpha}\left(\theta^{\alpha}\right)}\frac{\partial^{\alpha} f_1}{\partial \theta^{\alpha}}\right.\notag\\
	&-\frac{1}{r^{2\alpha}\sin_{\alpha}^{2\alpha}\left(\theta^{\alpha}\right)}\frac{\partial^{2\alpha} f_1}{\partial \psi^{2\alpha}}+\frac{2}{r^{\alpha}}\frac{\partial^{2\alpha} f_2}{\partial r^{\alpha} \partial \theta^{\alpha}}+\frac{2\cos_{\alpha}\left(\theta^{\alpha}\right)}{r^{\alpha}\sin_{\alpha}\left(\theta^{\alpha}\right)}\frac{\partial f_2}{\partial r^{\alpha}}\notag\\
	&\left.+\frac{2}{r^{\alpha}\sin_{\alpha}\left(\theta^{\alpha}\right)}\frac{\partial^{2\alpha} f_3}{\partial r^{\alpha} \partial \psi^{\alpha}}\right]\,\mathbf{e}_{r}^{\alpha}\notag\\
	&+\left[-\frac{\partial^{2\alpha} f_2}{\partial r^{2\alpha}}-\frac{2}{r^{\alpha}}\frac{\partial^{\alpha} f_2}{\partial r^{\alpha}}-\frac{1}{r^{2\alpha}\sin_{\alpha}^{2\alpha}\left(\theta^{\alpha}\right)}f_2+\frac{1}{r^{2\alpha}}\frac{\partial^{2\alpha} f_2}{\partial \theta^{2\alpha}}\right.\notag\\
	&+\frac{\cos_{\alpha}\left(\theta^{\alpha}\right)}{r^{2\alpha}\sin_{\alpha}\left(\theta^{\alpha}\right)}\frac{\partial^{\alpha} f_2}{\partial \theta^{\alpha}}-\frac{1}{r^{2\alpha}\sin_{\alpha}^{2\alpha}\left(\theta^{\alpha}\right)}\frac{\partial^{2\alpha} f_2}{\partial \psi^{2\alpha}}+\frac{2}{r^{2\alpha}}\frac{\partial^{\alpha} f_1}{\partial \theta^{\alpha}}+\frac{2}{r^{\alpha}}\frac{\partial^{2\alpha} f_1}{\partial r^{\alpha} \partial \theta^{\alpha}}\notag\\
	&\left.+\frac{2}{r^{2\alpha}\sin_{\alpha}\left(\theta^{\alpha}\right)}\frac{\partial^{2\alpha} f_3}{\partial \theta^{\alpha} \partial \psi^{\alpha}}\right]\,\mathbf{e}_{\theta}^{\alpha}\notag\\
	&+\left[-\frac{\partial^{2\alpha} f_3}{\partial r^{2\alpha}}-\frac{2}{r^{\alpha}}\frac{\partial^{\alpha} f_3}{\partial r^{\alpha}}+\frac{1}{r^{2\alpha}\sin_{\alpha}^{2\alpha}\left(\theta^{\alpha}\right)}f_3-\frac{1}{r^{2\alpha}}\frac{\partial^{2\alpha} f_3}{\partial \theta^{2\alpha}}\right.\notag\\
	&-\frac{\cos_{\alpha}\left(\theta^{\alpha}\right)}{r^{2\alpha}\sin_{\alpha}\left(\theta^{\alpha}\right)}\frac{\partial^{\alpha} f_3}{\partial \theta^{\alpha}}+\frac{1}{r^{2\alpha}\sin_{\alpha}^{2\alpha}\left(\theta^{\alpha}\right)}\frac{\partial^{2\alpha} f_3}{\partial \psi^{2\alpha}}+\frac{2}{r^{2\alpha}\sin_{\alpha}\left(\theta^{\alpha}\right)}\frac{\partial^{\alpha} f_1}{\partial \psi^{\alpha}}\notag\\
	&\left.+\frac{2}{r^{\alpha}\sin_{\alpha}\left(\theta^{\alpha}\right)}\frac{\partial^{2\alpha} f_1}{\partial r^{\alpha} \partial \psi^{\alpha}}+\frac{2}{r^{2\alpha}\sin_{\alpha}\left(\theta^{\alpha}\right)}\frac{\partial^{2\alpha} f_2}{\partial\theta^{\alpha}\partial \psi^{\alpha}}\right]\,\mathbf{e}_{\psi}^{\alpha}.
\end{align}

\section{The Cantor-type cylindrical and spherical-coordinate methods to the local fractional quaternionic Helmholtz equation}

In this section we derive interesting formulas for the component equations of the local fractional quaternionic Helmholtz equation
\begin{equation}\label{quaternionic Helmholtz}
	\Delta_{\mathbb H}^{(\alpha)}[f]+\lambda^2 f=0,\quad \lambda\in\mathbb{C},
\end{equation}
on the Cantor sets by using the Cantor-type cylindrical and spherical coordinate methods.

We have
\begin{equation}
	-\left(D_{MT}^{(\alpha)}-\lambda{\mathcal I}\right)\left(D_{MT}^{(\alpha)}+\lambda{\mathcal I}\right)=\Delta_{\mathbb H}^{(\alpha)}+\lambda^2{\mathcal I},
\end{equation}
where $\mathcal I$ denotes the identity operator.

Seen this factorization, the null solutions of the operator $D_{MT}^{(\alpha)}+\lambda{\mathcal I}$, called the perturbed local fractional Moisil-Teodorescu operator, are special solutions of (\ref{quaternionic Helmholtz}).

On substituting Eq. (\ref{scalar Laplacian cylindrical}) and Eq. (\ref{vector Laplacian cylindrical}) into Eq. (\ref{quaternionic Helmholtz}), we get the local fractional quaternionic Helmholtz equation in the Cantor-type cylindrical coordinates.

\begin{align}
0&=\left[\frac{\partial^{2\alpha} f_0}{\partial r^{2\alpha}}+\frac{1}{r^{2\alpha}}\frac{\partial^{2\alpha} f_0}{\partial \theta^{2\alpha}}+\frac{1}{r^{\alpha}}\frac{\partial^{\alpha} f_0}{\partial r^{\alpha}}+\frac{\partial^{2\alpha} f_0}{\partial z^{2\alpha}}+\lambda^2 f_0\right]\notag\\
&+\left[\frac{\partial^{2\alpha} f_1}{\partial r^{2\alpha}}+\frac{1}{r^{\alpha}}\frac{\partial^{\alpha} f_1}{\partial r^{\alpha}}-\frac{f_1}{r^{2\alpha}}+\frac{1}{r^{2\alpha}}\frac{\partial^{2\alpha} f_1}{\partial \theta^{2\alpha}}-\frac{2}{r^{2\alpha}}\frac{\partial^{\alpha} f_2}{\partial \theta^{\alpha}}+\frac{\partial^{2\alpha} f_1}{\partial z^{2\alpha}}+\lambda^2 f_1\right]\,\mathbf{e}_{r}^{\alpha}\notag\\	&+\left[\frac{\partial^{2\alpha} f_2}{\partial r^{2\alpha}}+\frac{1}{r^{\alpha}}\frac{\partial^{\alpha} f_2}{\partial r^{\alpha}}-\frac{f_2}{r^{2\alpha}}+\frac{1}{r^{2\alpha}}\frac{\partial^{2\alpha} f_2}{\partial \theta^{2\alpha}}+\frac{2}{r^{2\alpha}}\frac{\partial^{\alpha} f_1}{\partial \theta^{\alpha}}+\frac{\partial^{2\alpha} f_2}{\partial z^{2\alpha}}+\lambda^2 f_2\right]\,\mathbf{e}_{\theta}^{\alpha}\notag\\
&+\left[\frac{\partial^{2\alpha} f_3}{\partial r^{2\alpha}}+\frac{1}{r^{\alpha}}\frac{\partial^{\alpha} f_3}{\partial r^{\alpha}}+\frac{1}{r^{2\alpha}}\frac{\partial^{2\alpha} f_3}{\partial \theta^{2\alpha}}+\frac{\partial^{2\alpha} f_3}{\partial z^{\alpha}}+\lambda^2 f_3\right]\,\mathbf{e}_{z}^{\alpha}.\label{quaternionic Helmholtz cylindrical}
\end{align}

Equating each component in (\ref{quaternionic Helmholtz cylindrical}) to zero we get
\begin{equation}\label{quaternionic Helmholtz cylindrical_2}
	\begin{cases}
		\displaystyle \Delta_0^{(\alpha)}[f_0]+\lambda^2f_0=0,\\
		\displaystyle \Delta_0^{(\alpha)}[f_1]-\frac{2}{r^{2\alpha}}\frac{\partial^{\alpha} f_2}{\partial \theta^{\alpha}}+\left(\lambda^2 -\frac{1}{r^{2\alpha}}\right)f_1=0,\\
		\displaystyle \Delta_0^{(\alpha)}[f_2]+\frac{2}{r^{2\alpha}}\frac{\partial^{\alpha} f_1}{\partial \theta^{\alpha}}+\left(\lambda^2 -\frac{1}{r^{2\alpha}}\right)f_2=0,\\
		\displaystyle \Delta_0^{(\alpha)}[f_3]+\lambda^2f_3=0.
	\end{cases}
\end{equation}
One can notice that the components of the local fractional quaternionic Helmholtz operator in Cantor-type cylindrical coordinates are different from the corresponding local fractional scalar Helmholtz operator (there are several extra terms).

In a similar way, substituting Eq. (\ref{scalar Laplacian spherical}) and Eq. (\ref{vector Laplacian spherical}) into Eq. (\ref{quaternionic Helmholtz}) we obtain the local fractional quaternionic Helmholtz equation in the Cantor-type spherical coordinates.
\begin{align}
0&=\left[\frac{\partial^{2\alpha} f_0}{\partial r^{2\alpha}}+\frac{2}{r^{\alpha}}\frac{\partial^{\alpha} f_0}{\partial r^{\alpha}}+\frac{1}{r^{2\alpha}}\frac{\partial^{2\alpha} f_0}{\partial \theta^{2\alpha}}+\frac{\cos_{\alpha}\left(\theta^{\alpha}\right)}{r^{2\alpha}\sin_{\alpha}\left(\theta^{\alpha}\right)}\frac{\partial^{\alpha} f_0}{\partial \theta^{\alpha}}\right.\notag\\
&\left.+\frac{1}{r^{2\alpha}\sin_{\alpha}^{2\alpha}\left(\theta^{\alpha}\right)}\frac{\partial^{2\alpha} f_0}{\partial \psi^{2\alpha}}+\lambda^2 f_0\right]\notag\\
&+\left[\frac{\partial^{2\alpha} f_1}{\partial r^{2\alpha}}+\frac{2}{r^{\alpha}}\frac{\partial^{\alpha} f_1}{\partial r^{\alpha}}-\frac{2}{r^{2\alpha}}f_1+\frac{1}{r^{2\alpha}}\frac{\partial^{2\alpha} f_1}{\partial \theta^{2\alpha}}+\frac{\cos_{\alpha}\left(\theta^{\alpha}\right)}{r^{2\alpha}\sin_{\alpha}\left(\theta^{\alpha}\right)}\frac{\partial^{\alpha} f_1}{\partial \theta^{\alpha}}\right.\notag\\
&+\frac{1}{r^{2\alpha}\sin_{\alpha}^{2\alpha}\left(\theta^{\alpha}\right)}\frac{\partial^{2\alpha} f_1}{\partial \psi^{2\alpha}}-\frac{2}{r^{2\alpha}}\frac{\partial^{\alpha} f_2}{\partial \theta^{\alpha}}-\frac{2\cos_{\alpha}\left(\theta^{\alpha}\right)}{r^{2\alpha}\sin_{\alpha}\left(\theta^{\alpha}\right)}f_2\notag\\
&\left.-\frac{2}{r^{2\alpha}\sin_{\alpha}\left(\theta^{\alpha}\right)}\frac{\partial^{\alpha} f_3}{\partial \psi^{\alpha}}+\lambda^2 f_1\right]\,\mathbf{e}_{r}^{\alpha}\notag\\
&+\left[\frac{\partial^{2\alpha} f_2}{\partial r^{2\alpha}}+\frac{2}{r^{\alpha}}\frac{\partial^{\alpha} f_2}{\partial r^{\alpha}}-\frac{1}{r^{2\alpha}\sin_{\alpha}^{2\alpha}\left(\theta^{\alpha}\right)}f_2+\frac{1}{r^{2\alpha}}\frac{\partial^{2\alpha} f_2}{\partial \theta^{2\alpha}}+\frac{\cos_{\alpha}\left(\theta^{\alpha}\right)}{r^{2\alpha}\sin_{\alpha}\left(\theta^{\alpha}\right)}\frac{\partial^{\alpha} f_2}{\partial \theta^{\alpha}}\right.\notag\\
&+\left.\frac{1}{r^{2\alpha}\sin_{\alpha}^{2\alpha}\left(\theta^{\alpha}\right)}\frac{\partial^{2\alpha} f_2}{\partial \psi^{2\alpha}}+\frac{2}{r^{2\alpha}}\frac{\partial^{\alpha} f_1}{\partial \theta^{\alpha}}-\frac{2\cos_{\alpha}\left(\theta^{\alpha}\right)}{r^{2\alpha}\sin_{\alpha}^{2\alpha}\left(\theta^{\alpha}\right)}\frac{\partial^{\alpha} f_3}{\partial \psi^{\alpha}}+\lambda^2 f_2\right]\,\mathbf{e}_{\theta}^{\alpha}\notag\\
&+\left[\frac{\partial^{2\alpha} f_3}{\partial r^{2\alpha}}+\frac{2}{r^{\alpha}}\frac{\partial^{\alpha} f_3}{\partial r^{\alpha}}-\frac{1}{r^{2\alpha}\sin_{\alpha}^{2\alpha}\left(\theta^{\alpha}\right)}f_3+\frac{1}{r^{2\alpha}}\frac{\partial^{2\alpha} f_3}{\partial \theta^{2\alpha}}+\frac{\cos_{\alpha}\left(\theta^{\alpha}\right)}{r^{2\alpha}\sin_{\alpha}\left(\theta^{\alpha}\right)}\frac{\partial^{\alpha} f_3}{\partial \theta^{\alpha}}\right.\notag\\
&+\frac{1}{r^{2\alpha}\sin_{\alpha}^{2\alpha}\left(\theta^{\alpha}\right)}\frac{\partial^{2\alpha} f_3}{\partial \psi^{2\alpha}}+\frac{2}{r^{2\alpha}\sin_{\alpha}\left(\theta^{\alpha}\right)}\frac{\partial^{\alpha} f_1}{\partial \psi^{\alpha}}\notag\\
&\left.+\frac{2\cos_{\alpha}\left(\theta^{\alpha}\right)}{r^{2\alpha}\sin_{\alpha}^{2\alpha}\left(\theta^{\alpha}\right)}\frac{\partial^{\alpha} f_2}{\partial \psi^{\alpha}}+\lambda^2 f_3\right]\,\mathbf{e}_{\psi}^{\alpha}.\label{quaternionic Helmholtz spherical}
\end{align}

Now equate each component in (\ref{quaternionic Helmholtz spherical}) to zero we get
\begin{equation}\label{quaternionic Helmholtz spherical_2}
	\begin{cases}
		\displaystyle \Delta_0^{(\alpha)}[f_0]+\lambda^2f_0=0,\\
		\displaystyle \Delta_0^{(\alpha)}[f_1]-\frac{2}{r^{2\alpha}}\frac{\partial^{\alpha} f_2}{\partial \theta^{\alpha}}-\frac{2\cos_{\alpha}\left(\theta^{\alpha}\right)}{r^{2\alpha}\sin_{\alpha}\left(\theta^{\alpha}\right)}f_2-\frac{2}{r^{2\alpha}\sin_{\alpha}\left(\theta^{\alpha}\right)}\frac{\partial^{\alpha} f_3}{\partial \psi^{\alpha}}\\
		\displaystyle\qquad\qquad+\left(\lambda^2-\frac{2}{r^{2\alpha}}\right)f_1=0,\\
		\displaystyle \Delta_0^{(\alpha)}[f_2]+\frac{2}{r^{2\alpha}}\frac{\partial^{\alpha} f_1}{\partial \theta^{\alpha}}-\frac{2\cos_{\alpha}\left(\theta^{\alpha}\right)}{r^{2\alpha}\sin_{\alpha}^{2\alpha}\left(\theta^{\alpha}\right)}\frac{\partial^{\alpha} f_3}{\partial \psi^{\alpha}}\\
		\displaystyle\qquad\qquad+\left(\lambda^2 -\frac{1}{r^{2\alpha}\sin_{\alpha}^{2\alpha}\left(\theta^{\alpha}\right)}\right)f_2=0,\\
		\displaystyle \Delta_0^{(\alpha)}[f_3]+\frac{2}{r^{2\alpha}\sin_{\alpha}\left(\theta^{\alpha}\right)}\frac{\partial^{\alpha} f_1}{\partial \psi^{\alpha}}+\frac{2\cos_{\alpha}\left(\theta^{\alpha}\right)}{r^{2\alpha}\sin_{\alpha}^{2\alpha}\left(\theta^{\alpha}\right)}\frac{\partial^{\alpha} f_2}{\partial \psi^{\alpha}}\\
		\displaystyle\qquad\qquad+\left(\lambda^2 -\frac{1}{r^{2\alpha}\sin_{\alpha}^{2\alpha}\left(\theta^{\alpha}\right)}\right)f_3=0.
	\end{cases}
\end{equation}

Again, notice that the components of the local fractional quaternionic Helm\-holtz operator in Cantor-type spherical coordinates are different from the corresponding local fractional scalar Helmholtz operator.

Observe that the first component equation in (\ref{quaternionic Helmholtz cylindrical}) shows that the Helm\-holtz equation in Cantor-type cylindrical coordinates obtained in \cite{Hao} has exactly the form of the  scalar part of the local fractional quaternionic Helmholtz equation in the Cantor-type cylindrical coordinates. 

Similarly, first component equation in (\ref{quaternionic Helmholtz spherical}) reflects that the homogeneous version of the Helmholtz equation in Cantor-type spherical coordinates obtained in \cite{Rahmat} is identical with to  the scalar part of the local fractional quaternionic Helmholtz equation in the Cantor-type spherical coordinates.

\section*{Concluding Remarks}

In the present work, we have derived the local fractional quaternionic Moisil-Teodorescu and Helmholtz operators on the Cantor sets by using the Cantor-type cylindrical and spherical coordinates. The mathematical model developed in the present work emerges as a good generalization of the standard local fractional vector calculus in Cantor-type cylindrical and spherical coordinates.

\section*{Acknowledgements}
The authors were partially supported by Instituto Polit\'ecnico Nacional in the framework of SIP programs and by Fundaci\'{o}n Universidad de las Am\'{e}ricas Puebla, respectively.

\end{document}